\documentclass[12pt,a4paper]{amsart}
\usepackage[top=1.15in, bottom=1.15in, left=1.17in, right=1.17in]{geometry}
\usepackage{amssymb}
\usepackage{pdflscape}
\usepackage{amscd}
\usepackage{fancyhdr}

\newtheorem{theorem}{Theorem}[section]

\newtheorem{lemma}[theorem]{Lemma}
\newtheorem{proposition}[theorem]{Proposition}

\newtheorem{example}[theorem]{Example}
\newtheorem{definition}[theorem]{Definition}

\newtheorem{question}[theorem]{Question}

\title{Dimension expanders via quiver representations}
\author{Markus Reineke}
\begin{document}
\begin{abstract} We relate the notion of dimension expanders to quiver representations and their general subrepresentations, and use this relation to establish sharp existence results.
\end{abstract}
\maketitle
\parindent0pt

\section{Introduction}\label{introduction} Let $F$ be a field, and let $\varepsilon>0$ be a real number. An $\varepsilon$-expander is a tuple $(V,T_1,\ldots,T_k)$, consisting of a finite-dimensional $F$-vector space $V$, together with linear operators $T_1,\ldots,T_k$ on $V$, such that, for all subspaces $U\subset V$ of dimension $\dim U\leq \frac{1}{2}\dim V$, we have
$$\dim(U+\sum_{i=1}^kT_i(U))\geq(1+\varepsilon)\dim U.$$
This is a linear algebra analogue of the notion of expander graph \cite{L}. It is proven in \cite{LZ} for fields of characteristic zero, and in \cite{B,BY} for finite fields, that there exist $k$ and a {\it fixed} $\varepsilon>0$ such that $\varepsilon$-expanders $(V,T_1,\ldots,T_k)$ exist for {\it all} dimensions of the $F$-vector space $V$.\\[1ex]
In the present article, we sharpen this existence result and determine the optimal expansion coefficient $\varepsilon$ for $F$ an algebraically closed field.

\begin{theorem}\label{thmmain} Let $F$ be an algebraically closed field, let $k\geq 2$, and define $$\varepsilon_k=(k+1-\sqrt{k^2-2k+5})/2.$$ Then there exist $\varepsilon$-expanders $(V,T_1,\ldots,T_k)$ in all dimensions of $V$ if and only if $\varepsilon\leq\varepsilon_k$.
\end{theorem}

We will derive this result from a description of dimension vectors of subrepresentations of general representations of generalized Kronecker quivers; in particular, our proof will be non-constructive. However, this technique allows us to also cover the case of unbalanced dimension expanders of \cite{GRX}, and to speculate on potential generalizations of the notion of dimension expanders to arbitrary quivers.\\[1ex]
We will review the necessary quiver techniques in Section \ref{quivers}. In Section \ref{cdx}, we describe the dimension vectors of general representations of generalized Kronecker quivers. The applications to dimension expanders are derived in Section \ref{expanders}. Finally, in Section \ref{further}, we discuss potential generalizations to arbitrary quivers.\\[2ex]
{\bf Acknowledgments:} The author is grateful to M. Bertozzi and K. Martinez for many helpful comments on this manuscript. This work is supported by the DFG SFB-TRR 191 “Symplectic structures in geometry, algebra and dynamics”.

\section{Recollections on quiver representations}\label{quivers}
From now on, let $F$ be an algebraically closed field. 
Let $Q$ be a finite quiver with set of vertices $Q_0$ and arrows written $\alpha:i\rightarrow j$ for $i,j\in Q_0$, which we assume to be acyclic. We define the Euler form of $Q$ on $\mathbb{Z}Q_0$ by
$$\langle{\bf d},{\bf e}\rangle=\sum_{i\in Q_0}d_ie_i-\sum_{\alpha:i\rightarrow j}d_ie_j$$
for ${\bf d}=(d_i)_i$ and ${\bf e}=(e_i)_i$ in $\mathbb{Z}Q_0$. We consider the category ${\rm rep}_FQ$ of finite dimensional $F$-representations of $Q$, which is an abelian $F$-linear hereditary finite length category. Its Grothendieck group identifies with $\mathbb{Z}Q_0$ by associating to a representation $V$ its dimension vector ${\rm\bf dim}\, V$, and its homological Euler form is given by the Euler form:
$$\dim{\rm Hom}(V,W)-\dim{\rm Ext}^1(V,W)=\langle {\rm\bf dim}\, V,{\rm\bf dim}\, W\rangle$$
for all representations $V$ and $W$.\\[1ex]
For ${\bf d}\in\mathbb{N}Q_0$, we  fix $F$-vector spaces $V_i$ of dimension $d_i$ for all $i\in Q_0$. We define the representation space
$$R_{\bf d}(Q)=\bigoplus_{\alpha:i\rightarrow j}{\rm Hom}_F(V_i,V_j),$$
whose points $(f_\alpha)_\alpha$ we identify with the corresponding representation of $Q$ on the vector spaces $V_i$. On the $F$-vector space $R_{\bf d}(Q)$, the reductive linear algebraic group
$$G_{\bf d}=\prod_{i\in Q_0}{\rm GL}(V_i)$$
acts linearly via
$$(g_i)_i(f_\alpha)_\alpha=(g_jf_\alpha g_i^{-1})_{\alpha:i\rightarrow j}$$
such that the $G_{\bf d}$-orbits $\mathcal{O}_V$ in $R_{\bf d}(Q)$ naturally correspond to the isomorphism classes $[V]$ of $F$-representations $V$ of $Q$ of dimension vector ${\bf d}$.\\[1ex]
For ${\bf e}\leq{\bf d}$ componentwise, the subset of $R_{\bf d}(Q)$ of all representations $V$ admitting a subrepresentation of dimension vector ${\bf e}$ is Zariski-closed. Therefore, almost all representations in $R_{\bf d}(Q)$ (that is, those in a Zariski-dense subset) admit a subrepresentation of dimension vector ${\bf e}$ if and only if all representations in $R_{\bf d}(Q)$ do so. In this case, we write ${\bf e}\hookrightarrow{\bf d}$. There is a recursive numerical criterion for this notion due to Schofield in characteristic zero, generalized to positive characteristic by Crawley-Boevey:

\begin{theorem}[\cite{S,CB}]\label{generalsub} We have ${\bf e}\hookrightarrow{\bf d}$ if and only if $\langle{\bf e}',{\bf d}-{\bf e}\rangle\geq 0$ for all ${\bf e}'\hookrightarrow{\bf e}$.
\end{theorem}

\section{General subrepresentations of representation of generalized Kronecker quivers}\label{cdx}

Our main result in this section, which will directly apply to dimension expanders, is a non-recursive description of the relation ${\bf e}\hookrightarrow{\bf d}$ for generalized Kronecker quivers $K(m)$, given by two vertices $1$, $2$, and $m\geq 2$ arrows from $1$ to $2$. We prepare this description by some preliminary results.

\begin{lemma}\label{duality} For dimension vectors of $K(m)$, we have $(e_1,e_2)\hookrightarrow(d_1,d_2)$ if and only if $(d_2-e_2,d_1-e_1)\hookrightarrow(d_2,d_1)$.
\end{lemma}

\proof For a representation $V$ of $K(m)$ given by an $m$-tuple $f_1,\ldots,f_m:V_1\rightarrow V_2$ of linear maps, we denote by $V^*$ the representation $f_1^*,\ldots,f_m^*:V_2^*\rightarrow V_1^*$. This obviously defines a duality on ${\rm rep}_FK(m)$. Assume that a general representation $V$ of dimension vector ${\bf d}$ admits a subrepresentation $U$ of dimension vector ${\bf e}$. Then, dually, a general representation $V^*$ of dimension vector $(d_2,d_1)$ admits a factor representation of dimension vector $(e_2,e_1)$, whose kernel is a subrepresentation of dimension vector $(d_2-e_2,d_1-e_1)$. This finishes the proof.\\[1ex]
We now extend the Euler form $\langle\_,\_\rangle$ of $Q=K(m)$ to $\mathbb{R}Q_0$. We fix $0\not={\bf d}=(d_1,d_2)\in\mathbb{N}Q_0$ such that $0\geq\langle{\bf d},{\bf d}\rangle=d_1^2+d_2^2-md_1d_2$. In particular, $d_1,d_2\geq 1$, and $$(m-\sqrt{m^2-4})/2\leq {d_2}/{d_1}\leq ({m+\sqrt{m^2-4}})/2=:\beta.$$
For fixed $x\in[0,d_1]$, we consider the function
$$q_x(y)=\langle(x,y),(d_1-x,d_2-y)\rangle$$
on $[0,d_2]$, and denote by $c_{\bf d}$  the smaller of its two zeroes. The explicit form
$$c_{\bf d}(x)=\left(mx+d_2-\sqrt{(mx-d_2)^2+4x(d_1-x)}\right)/2$$
shows existence. In particular, $c_{\bf d}(x)\leq({mx+d_2})/{2}$, and ${mx+d_2}-c_{\bf d}(x)$ is the larger zero of $q_x$. We have $q_x(y)\geq 0$ for $c_{\bf d}(x)\leq y\leq mx+d_2-c_{\bf d}(x)$, and $q_x(y)\leq 0$ otherwise. We have the following estimate:

\begin{lemma}\label{estimate} If $\langle{\bf d},{\bf d}\rangle\leq 0$, we have ${d_2}/{d_1}\cdot x\leq c_{\bf d}(x)\leq \min(mx, d_2)$.
\end{lemma}

{\bf Proof:} We have $$q_x({d_2}/{d_1}\cdot x)=\langle{x}/{d_1}\cdot{\bf d},{\bf d}-{x}/{d_1}\cdot{\bf d}\rangle={x(d_1-x)}/{d_1^2}\cdot\langle{\bf d},{\bf d}\rangle\leq 0$$ by assumption. Thus the first inequality follows, since $c_{\bf d}(x)\leq({mx+d_2})/{2}$,  once we know  that ${d_2}/{d_1}\cdot x\leq({mx+d_2})/{2}$. If ${d_2}/{d_1}\leq{m}/{2}$, this holds trivially. Otherwise, we use $d_2\leq md_1$ to estimate
$$({d_2}/{d_1}-{m}/{2})x\leq({d_2}/{d_1}-{m}/{2})d_1={d_2}/{2}+({d_2-md_1})/{2}\leq{d_2}/{2},$$
and again the desired estimate follows.\\[1ex]
For the second inequality, we calculate
$$q_x(mx)=\langle x\cdot(1,m),{\bf d}-x\cdot(1,m)\rangle=x(d_1-x)\geq 0,$$
thus
$$c_{\bf d}(x)\leq mx\leq mx+d_2-c_{\bf d}(x),$$
which finishes the proof.

\begin{lemma}\label{preprojective} If $d_2>\beta d_1$ and ${\bf e}\hookrightarrow{\bf d}$, then $e_2>\beta e_1$.
\end{lemma}

\proof ${\bf e}\hookrightarrow{\bf d}$ implies $\langle{\bf e},{\bf d}-{\bf e}\rangle\geq 0$ by Schofield's criterion, thus $e_2\geq c_{\bf d}(e_1)$ by definition of $c_{\bf d}$. It thus suffices to prove that $c_{\bf d}(x)>\beta x$ provided $d_2>\beta d_1$. Since $\beta^2-m\beta+1=0$, we have $$q_x(\beta x)=x(d_1-(m-\beta)d_2)<xd_1(1-m\beta+\beta^2)=0,$$
from which we can conclude $c_{\bf d}(x)>\beta x$ provided $\beta x\leq mx+d_2-c_{\bf d}(x)$. But $$\beta x<mx\leq mx+d_2-c_{\bf d}(x)$$
since $c_{\bf d}(x)\leq d_2$.\\[1ex]
We can now derive the main result of this section:

\begin{proposition}\label{propmain} If $Q=K(m)$ is the $m$-arrow Kronecker quiver $\bullet\stackrel{(m)}{\Rightarrow}\bullet$ and $\langle{\bf d},{\bf d}\rangle\leq 0$, then for ${\bf e}\leq{\bf d}$ the following are equivalent:
\begin{enumerate}
\item ${\bf e}\hookrightarrow{\bf d}$,
\item $\langle{\bf e},{\bf d}-{\bf e}\rangle\geq 0$,
\item $e_2\geq c_{\bf d}(e_1)$.
\end{enumerate}
\end{proposition}

\proof Without loss of generality, we assume $d_1\leq d_2$ using the duality Lemma \ref{duality}. Obviously (1) implies (2) implies (3): if ${\bf e}\hookrightarrow{\bf d}$ then, by Schofield's criterion applied to ${\bf e}'={\bf e}$, we have $\langle{\bf e},{\bf d}-{\bf e}\rangle\geq 0$, and by definition of the function $c_{\bf d}$, this implies $e_2\geq c_{\bf d}(e_1)$. Conversely, assume that this inequality holds, and let ${\bf e'}\hookrightarrow{\bf e}$. We want to prove that $\langle{\bf e}',{\bf d}-{\bf e}\rangle\geq 0$; then ${\bf e}\hookrightarrow{\bf d}$ follows from Schofield's criterion. We first assume $\langle{\bf e},{\bf e}\rangle\geq 1$, thus $e_2>\beta e_1$ or $e_2<(m-\beta)e_1$. Since $e_2\geq c_{\bf d}(e_1)\geq {d_2}/{d_1}\cdot e_1\geq e_1$ by assumption and Lemma \ref{estimate}, we have $e_2>\beta e_1$ since $\beta\geq 1$. By Lemma \ref{preprojective}, we find $e_2'>\beta e_1'$, and thus
$$\langle{\bf e'},{\bf d}-{\bf e}\rangle=e_1'(d_1-e_1-m(d_2-e_2))+e_2'(d_2-e_2)>e_1'(d_1-e_1-(m-\beta)(d_2-e_2)).$$
Since $d_2\leq\beta d_1$, we have $d_2-e_2<\beta(d_1-e_1)$, thus $d_1-e_1>(m-\beta)(d_2-e_2)$, proving the claim. Now we assume $\langle{\bf e},{\bf e}\rangle\leq 0$. By Lemma \ref{estimate}, we have $$e_2'\geq c_{\bf e}(e_1')\geq{e_2}/{e_1}\cdot e_1',$$
and thus
$$\langle{\bf e'},{\bf d}-{\bf e}\rangle=e_1'(d_1-e_1-m(d_2-e_2))+e_2'(d_2-e_2)\geq$$
$$\geq{e_1'}/{e_1}\cdot\left(e_1(d_1-e_1)-me_1(d_2-e_2)+e_2(d_2-e_2)\right)={e_1'}/{e_1}\cdot\langle{\bf e},{\bf d}-{\bf e}\rangle\geq 0,$$
again proving the claim. 

\section{Application to dimension expanders}\label{expanders}

We generalize the definition of dimension expanders of Section \ref{introduction} to a notion of expander representation. Proposition \ref{propmain} then almost immediately yields sharp existence results.

\begin{definition} Let $0<\delta<1$ and $\varepsilon>0$, and let $V$ and $W$ be non-zero finite-dimensional $F$-vector spaces. We call a representation $f_1,\ldots,f_m:V\rightarrow W$ of $K(m)$ a $(\delta,\varepsilon)$-expander representation if for all subspaces $0\not=U\subset V$ such that $\frac{\dim U}{\dim V}\leq\delta$, we have
$$\dim\sum_{k=1}^mf_k(U)\geq (1+\varepsilon)\cdot\frac{\dim W}{\dim V}\cdot\dim U.$$
\end{definition}

The following lemma translates existence of expander representations to properties of dimension vectors of subrepresentations of general representations:

\begin{lemma} For fixed integers $m,d_1,d_2 \geq 1$ and real numbers $0<\delta<1$, $\varepsilon>0$, there exists a $(\delta,\varepsilon)$-expander representation of $K(m)$ of dimension vector $(d_1,d_2)$ if and only if for all $(e_1,e_2)\hookrightarrow(d_1,d_2)$ such that $e_1\leq\delta\cdot d_1$, we have
$$e_2\geq(1+\varepsilon)\cdot\frac{d_2}{d_1}\cdot e_1.$$
\end{lemma}

\proof Assume there exists such an expander representation $M$ given by $f_1,\ldots,f_m:V\rightarrow W$, and assume that $(e_1,e_2)\hookrightarrow(d_1,d_2)$. Then in particular $M$ admits a subrepresentation of dimension vector $(e_1,e_2)$, that is, there exists a subspace $U\subset V$ of dimension $e_1$ such that $\sum_kf_k(U)$ is of dimension at most $e_2$. On the other hand, $\sum_kf_k(U)$ is at least of dimension $(1+\varepsilon)\cdot\frac{d_2}{d_1}\cdot\dim U$. The claimed inequality for $e_2$ follows. Conversely, assume that the numerical condition is satisfied. Then the set $S_{(e_1,e_2)}\subset R_{(d_1,d_2)}(K(m))$ of representations admitting a subrepresentation of dimension vector $(e_1,e_2)$ is a proper Zariski-closed subset whenever $e_2<(1+\varepsilon)\cdot\frac{d_2}{d_1}\cdot e_1$. Thus the union of all these finitely many proper closed subsets is again a proper subset, and any representation in its complement is a $(\delta,\varepsilon)$-expander representation by definition.\\[1ex]
This allows us to establish the following sharp existence result:

\begin{theorem} Fix an integer $m\geq 1$, real numbers $0<\delta<1$ and $\varepsilon>0$, and a rational $\alpha$ such that $\alpha^2-m\alpha+1<0$ and $m\delta+\alpha-2\alpha\delta>0$. Define 

$$\varepsilon_m(\alpha,\delta)=\frac{m\delta+\alpha-2\alpha\delta-\sqrt{(m\delta-\alpha)^2+4\delta(1-\delta)}}{2\alpha\delta}>0.$$

Then there exist $(\delta,\varepsilon)$-expander representations of $K(m)$ for all dimension vectors $(d_1,d_2)$ such that $d_2/d_1=\alpha$ if and only if $\varepsilon\leq\varepsilon_m(\alpha,\delta)$.
\end{theorem}

\proof The assumptions on $\alpha$ ensure that $\varepsilon_m(\alpha,\delta)>0$ by a straightforward calculation. We consider dimension vectors ${\bf d}$ such that $d_2/d_1=\alpha$; in particular $\langle{\bf d},{\bf d}\rangle< 0$. By the previous lemma and Proposition \ref{propmain}, we have:\\[1ex]
There exists a $(\delta,\varepsilon)$-expander representation of $K(m)$ of dimension vector ${\bf d}$ if and only if $e_2\geq (1+\varepsilon)\alpha e_1$ for all $e_1\leq \delta d_1$ and all $e_2\geq c_{\bf d}(e_1)$, or, equivalently, if $\lceil c_{\bf d}(x)\rceil\geq(1+\varepsilon)\alpha x$ for all integral $x\leq\delta\cdot d_1$.\\[1ex]
This implies:\\[1ex]
There exist $(\delta,\varepsilon)$-expander representations of $K(m)$ for all dimension vectors $(d_1,d_2)$ such that $d_2/d_1=\alpha$ if and only if $\lceil c_{{\bf d}}(x)\rceil\geq(1+\varepsilon)\alpha x$ for all dimension vectors ${\bf d}=(d_1,d_2)$ such that $d_2/d_1=\alpha$ and all integral $x\leq \delta d_1$.\\[1ex]
The function $c_{\bf d}(x)$ is concave on the interval $[0,d_2]$ since, by a straightforward calculation, its second derivative equals
$$c_{\bf d}''(x)=\frac{2\langle{\bf d},{\bf d}\rangle}{\left({(mx-d_2)^2+4x(d_1-x)}\right)^{3/2}},$$
which is negative by assumption. Thus, in the interval $[0,\delta d_1]$, the fraction $c_{\bf d}(x)/x$ attains its minimum at $\delta d_1$. For $\rho\in[0,1]$, we have ${c_{\bf d}(\rho d_1)}/{(\alpha\rho d_1)}=1+\varepsilon_m(\alpha,\rho),$
thus in the interval $[0,\delta]$, the function $\varepsilon_m(\alpha,\rho)$ of $\rho$ attains its minimum at $\rho=\delta$. We thus find that the above existence condition is equivalent to 
$$\lceil (1+\varepsilon_m(\alpha,\rho))\alpha\rho d_1\rceil\geq(1+\varepsilon)\alpha\rho d_1$$
for all $d_1$ such that $\alpha d_1$ is integral and all $\rho\in[0,\delta]$ such that $\rho d_1$ is integral. This is clearly equivalent to $\varepsilon_m(\alpha,\rho)\geq\varepsilon$ for all $\rho\in[0,\delta]$, and this in turn to $\varepsilon_m(\alpha,\delta)\geq\varepsilon$. This finishes the proof.\\[1ex]
This result immediately implies Theorem \ref{thmmain} as the special case $m=k+1$, $\alpha=1$, $\delta=1/2$. Namely, in a general representation of $K(k+1)$ of dimension vector $(d,d)$, the map representing the first arrow is invertible, thus w.l.o.g.~the identity, and ${\rm id},T_1,\ldots,T_k:V\rightarrow V$ defines an expander representation if and only if $(V,T_1,\ldots,T_k)$ is a dimension expander; moreover, $\varepsilon_{k+1}(1,1/2)=\varepsilon_k$.

\section{Potential generalizations}\label{further}

We finish with a few remarks suggesting further directions.\\[1ex]
The characterization of dimension vectors ${\bf e}\hookrightarrow{\bf d}$ of subrepresentations of general representation by the single quadratic equation $\langle{\bf e},{\bf d}-{\bf e}\rangle\geq 0$ of Proposition \ref{propmain} is special to the quivers $K(m)$. Namely, we have the following

\begin{example} For the complete bipartite three-vertex quiver $\bullet\Rightarrow\bullet\Leftarrow\bullet$, the dimension vector ${\bf d}=(3,6,5)$ is a Schur root (even belonging to the fundamental domain), and ${\bf e}=(3,5,1)$ fulfills $\langle{\bf e},{\bf d}-{\bf e}\rangle\geq 0$ (even $>0$), but ${\bf e}\not\hookrightarrow{\bf d}$.\end{example}

It is natural to ask whether the explicit dimension expanders constructed in \cite{LZ} using representations of ${\rm SL}_2(\mathbb{Z})$, and in \cite{B} using monotone expanders, are already $\varepsilon_k$-expanders for the optimal expansion coefficients $\varepsilon_k$.\\[1ex]
In another direction, dimension expanders were used in \cite{Eck} to construct non-hyper\-finite families of representations of generalized Kronecker quivers, and it would be interesting to know whether the present methods yield new insights about such families.\\[1ex] 
Representations $f_1,\ldots,f_m:V_1\rightarrow V_2$ such that $\dim\sum_kf_k(U)>\frac{\dim V_2}{\dim V_1}\dim U$ for all proper non-zero subspaces $U\subset V_1$ are stable in the sense of Geometric Invariant Theory \cite{King}, thus the $(\delta,\varepsilon)$-expander property might be viewed as a quantitative form of stability. This point of view suggests a generalization to arbitrary quivers:

\begin{definition} Let $Q$ be a finite quiver, and let $\Theta\in(\mathbb{R}Q_0)^*$ be a stability function for $Q$. Let ${\bf d}\in\mathbb{N}Q_0$ be a dimension vector for $Q$ such that $\Theta({\bf d})=0$, and let $0<\delta<1$ and $\varepsilon>0$ be reals. We call a representation $V$ of $Q$ of dimension vector ${\bf\dim} V={\bf d}$ a $(\delta,\varepsilon)$-expander relative to $\Theta$ if for all subrepresentations $U\subset V$ such that $\dim U\leq\delta\cdot\dim V$, we have $\Theta({\rm\bf dim}\, U)\leq-\varepsilon\cdot\dim U$.
\end{definition}

Then it is natural to ask for uniform expansion:

\begin{question} For which $\Theta$, $\delta$ and $\varepsilon$ is it true that there exist $(\delta,\varepsilon)$-expander representations relative to $\Theta$ for all (resp. all sufficiently large) dimension vectors ${\bf d}$ such that $\Theta({\bf d})=0$?
\end{question}

\end{document}